\documentclass[11pt,a4paper,reqno]{amsart}
\usepackage[margin=1.18in]{geometry}
\usepackage{amsmath,amssymb,amsthm,mathtools,bm,mathrsfs}
\usepackage{float}
\usepackage[colorlinks=true,linkcolor=black,citecolor=black,urlcolor=black]{hyperref}
\hypersetup{pdftitle={The distribution of greatest common divisors of binomial coefficients and their upper indices},pdfauthor={Zhoulin Cui and Zhishan Yang}}

\usepackage[no-math]{fontspec}
\allowdisplaybreaks[2]
\numberwithin{equation}{section}

\newtheoremstyle{paperplain}{6pt}{6pt}{\itshape}{0pt}{\bfseries}{.}{0.5em}{}
\theoremstyle{paperplain}
\newtheorem{theorem}{Theorem}[section]
\newtheorem{lemma}[theorem]{Lemma}
\newtheorem{proposition}[theorem]{Proposition}
\newtheorem{corollary}[theorem]{Corollary}
\theoremstyle{remark}\newtheorem{remark}[theorem]{Remark}

\newcommand{\ind}{\mathbf 1}
\newcommand{\cc}{\mathfrak c}

\title[Greatest common divisors of binomial coefficients]{The distribution of greatest common divisors of binomial coefficients and their upper indices}
\author{Zhoulin Cui and Zhishan Yang}
\address[Zhoulin Cui]{School of Mathematics and Statistics, Qingdao University, Qingdao, Shandong 266071, P.R. China}
\email{\href{mailto:cuizl02@163.com}{cuizl02@163.com}}
\address[Zhishan Yang, corresponding author]{School of Mathematics and Statistics, Qingdao University, Qingdao, Shandong 266071, P.R. China}
\email{\href{mailto:zsyang@qdu.edu.cn}{zsyang@qdu.edu.cn}}
\date{}

\begin{document}

\begin{abstract}
For positive integers $n,d$, let $A_d(n)$ denote the number of integers $0\le k\le n$ such that $\gcd\bigl(n,\binom nk\bigr)=d$, and put $S_d(x)=\sum_{n\le x}A_d(n)$. We prove that, for each fixed integer $d\ge2$,
\[
 S_d(x)\sim\frac{d-1}{d}\,\cc\,\frac{x}{\log x}\qquad(x\to\infty),\qquad \cc=1.526453\ldots
\]
The key step is to extend the coprimality result of Ford and Konyagin to arbitrary fixed parameters $d\ge2$ and $1\le j<d$, obtaining
\[
 \#\Bigl\{m\le x:\gcd\Bigl(m,\binom{dm}{jm}\Bigr)=1\Bigr\}
 \sim\cc\frac x{\log x}.
\]
The zero set of the corresponding Kummer carry function is a finite union of intervals of total measure $1/2$. We combine this property with estimates for small prime factors, Selberg trigonometric approximation, and exponential sum estimates over primes to count simultaneous carry conditions. We further prove that, for every fixed integer $e\ge1$,
\[
 \#\Bigl\{m\le x:\gcd\Bigl(m,\binom{dm}{jm}\Bigr)=e\Bigr\}
 \sim\frac{\cc}{e}\frac x{\log x}.
\]
\end{abstract}

\maketitle

\section{Introduction and main results}\label{sec1}

The divisibility properties of the central binomial coefficient $\binom{2n}{n}$ exhibit a simple but striking contrast. The integrality of the Catalan number $\frac1{n+1}\binom{2n}{n}$ ensures that $n+1$ always divides $\binom{2n}{n}$, whereas divisibility by the adjacent integer $n$ depends on the arithmetic structure of $n$. Pomerance~\cite{bib14} studied such questions systematically. This leads naturally to the greatest common divisor of $n$ and $\binom{2n}{n}$: divisibility and coprimality describe the two extreme cases, while prescribing the greatest common divisor gives a finer counting problem.

The basic tools for these questions arise from the connection between binomial coefficients and digital representations. In 1852, Kummer~\cite{bib9} proved that $v_p\bigl(\binom{a+b}{a}\bigr)$ equals the number of carries when $a$ and $b$ are added in base $p$. Lucas's congruence theorem~\cite{bib10}, established in 1878, expresses a binomial coefficient modulo $p$ as a product of binomial coefficients formed from the corresponding base-$p$ digits. For arithmetic properties modulo prime powers and for extensions and applications of Lucas's theorem, see Granville~\cite{bib5} and Meštrović~\cite{bib12}, respectively. These results turn prime divisibility into restrictions on digits, making carries a basic link between binomial coefficients and counting problems.

Digital methods also yield divisibility statistics in Pascal's triangle. Fine~\cite{bib3} obtained a formula for the number of binomial coefficients in a row that are not divisible by a given prime, and Carlitz~\cite{bib1} investigated the corresponding counts for fixed prime powers. In a different direction, Joris, Oestreicher, and Steinig~\cite{bib8} studied the greatest common divisor of a collection of binomial coefficients with consecutive lower indices, while McTague~\cite{bib11} considered common factors when the lower indices lie in an arithmetic progression. These works describe the arithmetic of binomial coefficients through prime valuations or common factors of a collection of coefficients. The greatest common divisor of an individual coefficient and its upper index involves both the carry structure of the coefficient and the prime factorization of the upper index.

For central binomial coefficients, Erdős, Graham, Ruzsa, and Straus~\cite{bib2} studied divisibility by high powers of small primes and proved that, for any two given odd primes $p,q$, there are infinitely many $n$ for which $\binom{2n}{n}$ is coprime to $pq$. Granville and Ramaré~\cite{bib6} applied explicit exponential sum estimates to the squarefreeness of binomial coefficients, demonstrating the role of analytic methods in controlling divisibility conditions at several primes. Replacing the fixed integer $pq$ in the coprimality condition by the varying index $n$ requires simultaneous control of all prime factors of $n$, and hence a combination of digital restrictions with information on the distribution of those prime factors.

For the set of integers satisfying $n\mid\binom{2n}{n}$, Pomerance~\cite[Theorem~3 and its proof]{bib14} obtained the upper density bound $1-\log 2$, which was subsequently improved by Sanna~\cite[Theorem~1.2]{bib15}. In the coprime case, Sanna~\cite[Theorem~1.3]{bib15} proved that the number of $n\le x$ satisfying $\gcd\bigl(n,\binom{2n}{n}\bigr)=1$ is $O(x/\sqrt{\log x})$. Since every odd prime satisfies this condition, there is also a lower bound of $(1+o(1))x/\log x$. Ford and Konyagin~\cite[Theorems~1, 3]{bib4} later proved that, for each fixed positive integer $\ell$, the set of integers satisfying $n^\ell\mid\binom{2n}{n}$ has positive asymptotic density. They also determined the leading constant in the coprime count:
\begin{equation}
 \#\Bigl\{n\le x:\gcd\Bigl(n,\binom{2n}{n}\Bigr)=1\Bigr\}
 \sim\cc\,\frac{x}{\log x},\qquad \cc=1.526453\ldots.
 \label{eq1.1}
\end{equation}

Formula~\eqref{eq1.1} suggests two natural questions. How does the leading constant change when the central position is replaced by an arbitrary fixed ratio of lower to upper indices? How does the main term depend on the prescribed value when coprimality is replaced by a fixed greatest common divisor? We establish the corresponding asymptotic formulas and use them to count pairs of indices for which a binomial coefficient and its upper index have a given greatest common divisor. This requires treating carry conditions that vary with the index ratio and controlling the local restrictions imposed by prescribed common prime factors.

For a fixed positive integer $d$, we study the number of integer pairs $(n,k)$ satisfying
\[
 \gcd\Bigl(n,\binom nk\Bigr)=d,\qquad 1\le n\le x,\quad 0\le k\le n,
\]
To this end, write
\[
 A_d(n)=\#\Bigl\{0\le k\le n:\gcd\Bigl(n,\binom nk\Bigr)=d\Bigr\},\qquad
 S_d(x)=\sum_{n\le x}A_d(n).
\]
When $d=1$, only the endpoints $k=0,n$ satisfy the condition, so $A_1(n)=2$. When $d\ge2$, the count comes from the interior indices, and its average behaviour is related to coprimality properties of central binomial coefficients.

Our first main result is as follows.

\begin{theorem}\label{thm1.1}
For all $x\ge1$, we have $S_1(x)=2\lfloor x\rfloor$. For each fixed integer $d\ge2$,
\begin{equation}
 S_d(x)\sim\frac{d-1}{d}\,\cc\,\frac{x}{\log x}\qquad(x\to\infty).
 \label{eq1.2}
\end{equation}
\end{theorem}

The identity $k\binom nk=n\binom{n-1}{k-1}$ implies that, if $d\ge2$ and $\gcd\bigl(n,\binom nk\bigr)=d$, then $n=dm$ and $k=jm$ for some $1\le j<d$. Thus the problem reduces to counting binomial coefficients with a fixed ratio of indices. For fixed integers $d\ge2$, $1\le j<d$, and $e\ge1$, let
\[
 B_{d,j}(x;e)=\#\Bigl\{m\le x:\gcd\Bigl(m,\binom{dm}{jm}\Bigr)=e\Bigr\},
 \qquad B_{d,j}(x)=B_{d,j}(x;1).
\]
The following result covers both coprimality at a general fixed ratio and a prescribed greatest common divisor.

\begin{theorem}\label{thm1.2}
For fixed integers $d\ge2$, $1\le j<d$, and $e\ge1$, we have
\begin{equation}
 B_{d,j}(x;e)\sim\frac{\cc}{e}\frac x{\log x}\qquad(x\to\infty).
 \label{eq1.3}
\end{equation}
In particular, when $e=1$,
\begin{equation}
 B_{d,j}(x)\sim\cc\,\frac x{\log x}\qquad(x\to\infty).
 \label{eq1.4}
\end{equation}
\end{theorem}

Taking $(d,j,e)=(2,1,1)$ recovers the formula of Ford and Konyagin in~\eqref{eq1.1}. For $(d,j)=(2,1)$ and any fixed $e\ge1$, the theorem gives
\begin{equation}
 \#\Bigl\{n\le x:\gcd\Bigl(n,\binom{2n}{n}\Bigr)=e\Bigr\}
 \sim\frac{\cc}{e}\frac{x}{\log x}\qquad(x\to\infty).
 \label{eq1.5}
\end{equation}

The two theorems exhibit different types of dependence on the parameters. For each fixed pair $(d,j)$, the ratio $j/d$ changes the shape of the zero set of the carry function, but the leading constant in the coprime count remains the same constant $\cc$. Prescribing the greatest common divisor $e$ introduces a factor of $1/e$ in the main term. The key is that the general carry zero set is still a finite union of intervals of measure $1/2$, and simultaneous carry conditions retain the main-term weights from the central binomial case. For the original count of integer pairs, the proportionality forced by the fixed greatest common divisor yields the additional factor $(d-1)/d$. Moreover, when the count is split according to whether $n/d$ is prime or composite, the composite part also contributes a nonzero main term of the same order as the prime part; see Corollary~\ref{cor6.3}.

We follow the analytic framework of Ford and Konyagin~\cite{bib4}. Section~\ref{sec2} gives the counting reduction and introduces the carry function. Section~\ref{sec3} establishes valuation estimates at fixed primes and a bound for the contribution of small prime factors. Section~\ref{sec4} uses trigonometric approximation and exponential sum estimates over primes to count simultaneous carry conditions. In Section~\ref{sec5}, we compute the main term for each number of prime factors and sum the resulting contributions to obtain the coprime asymptotic~\eqref{eq1.4}. Section~\ref{sec6} completes the proofs of the two main theorems and gives applications. For general $e$, after the substitution $m=et$, the resulting condition can differ from the coprimality condition with parameters $(de,je)$ only at primes $p\mid e$; Proposition~\ref{prop6.2} bounds this difference with a power saving.

\medskip

\noindent\textbf{Notation.}
The letters $p,q$, with or without subscripts, denote primes. We write $\{t\}$ for the fractional part of a real number $t$,
and $v_p(n)$ for the $p$-adic valuation of a positive integer $n$. The implied constants in $\ll$ and $O(\cdot)$
may depend on the parameters indicated by subscripts.

\section{A counting reduction and the carry function}\label{sec2}

We first reduce the count of pairs $(n,k)$ to a count with a fixed ratio of indices, and then express the valuations at all primes in terms of a carry function.

\begin{lemma}\label{lem2.1}
For $d\ge2$ and $x\ge1$,
\begin{equation}
 S_d(x)=\sum_{j=1}^{d-1}T_{d,j}(x/d),\qquad
 T_{d,j}(x)=\#\Bigl\{m\le x:\gcd\Bigl(dm,\binom{dm}{jm}\Bigr)=d\Bigr\}.
 \label{eq2.1}
\end{equation}
\end{lemma}

\begin{proof}
Suppose that $\gcd\bigl(n,\binom nk\bigr)=d$. Since $d\mid n$, we may write $n=dm$ and $\binom nk=dB$, where $\gcd(m,B)=1$. As $d\ge2$, the endpoints $k=0,n$ do not satisfy the condition. For $1\le k<n$, the identity $k\binom nk=n\binom{n-1}{k-1}$ gives $kB=m\binom{n-1}{k-1}$, and hence $m\mid k$. Thus $k=jm$ with $1\le j<d$. Conversely, each pair $(m,j)$ satisfying $\gcd\bigl(dm,\binom{dm}{jm}\bigr)=d$ gives a pair $(n,k)$ of the required type. Summing over $j$ yields~\eqref{eq2.1}.
\end{proof}

Fix $d\ge2$ and $1\le j<d$, and abbreviate $C(m)=\binom{dm}{jm}$. For central binomial coefficients, the no-carry condition corresponds to the half-interval $[0,1/2)$. In the general case, define
\[
 F(t)=\lfloor dt\rfloor-\lfloor jt\rfloor-\lfloor(d-j)t\rfloor,
 \qquad Z=\{t\in[0,1):F(t)=0\}.
\]

\begin{lemma}\label{lem2.2}
The function $F$ is $1$-periodic and takes only the values $0,1$. Its zero set $Z$ is a union of at most $2d+1$ intervals and has Lebesgue measure $1/2$. For every prime $p$ and positive integer $m$,
\begin{equation}
 v_p(C(m))=\sum_{r\ge1}F\Bigl(\Bigl\{\frac m{p^r}\Bigr\}\Bigr).
 \label{eq2.2}
\end{equation}
Moreover, $F$ is identically $0$ on $[0,1/d)$ and identically $1$ on $I_d=\bigl(1/d,1/(d-1)\bigr)$.
\end{lemma}

\begin{proof}
The range of $F$ follows from $\lfloor a+b\rfloor-\lfloor a\rfloor-\lfloor b\rfloor\in\{0,1\}$, and its $1$-periodicity follows from $d=j+(d-j)$. Its discontinuities lie in $\frac1d\mathbb Z\cup\frac1j\mathbb Z\cup\frac1{d-j}\mathbb Z$, so $Z$ is a union of at most $2d+1$ intervals. If none of $dt,jt,(d-j)t$ is an integer, the identity $\lfloor-y\rfloor=-\lfloor y\rfloor-1$ gives $F(1-t)=1-F(t)$. Apart from finitely many points, the map $t\mapsto1-t$ interchanges the zero set of $F$ and its complement. Hence $|Z|=1/2$.

For $0\le t<1/d$, all three floor terms defining $F$ are $0$. For $t\in I_d$, we have $1<dt<d/(d-1)\le2$, while $jt<1$ and $(d-j)t<1$, so $F(t)=1$. Finally, Legendre's formula gives
\[
 v_p(C(m))=\sum_{r\ge1}\Bigl(\Bigl\lfloor\frac{dm}{p^r}\Bigr\rfloor-\Bigl\lfloor\frac{jm}{p^r}\Bigr\rfloor-\Bigl\lfloor\frac{(d-j)m}{p^r}\Bigr\rfloor\Bigr),
\]
Each summand is $F(\{m/p^r\})$, which proves~\eqref{eq2.2}.
\end{proof}

\section{Local valuations and small prime factors}\label{sec3}

A carry occurs throughout the interval $I_d$, whose length we denote by
\[
 \ell_d=|I_d|=\frac1{d(d-1)}.
\]
If $v_p(C(m))$ is small, a suitable grouping of the base-$p$ digits of $m$ restricts the number of occurrences of certain digits.
We first use this restriction to obtain a counting estimate with a power saving for each fixed prime $p$.

\begin{lemma}\label{lem3.1}
For each fixed prime $p$, there is a constant $\eta_p>0$, depending only on $d,p$, such that, for every fixed nonnegative integer $A$,
\begin{equation}
 \#\{m\le x:v_p(C(m))\le A\}
 \ll_{d,p,A}x^{1-\eta_p}(1+\log x)^A.
 \label{eq3.1}
\end{equation}
\end{lemma}

\begin{proof}
Choose a positive integer $L_d$, depending only on $d$, such that $2^{L_d}\ge4/\ell_d$, and put $M=p^{L_d}$. This amounts to combining each block of $L_d$ consecutive base-$p$ digits into a single digit. Select the digit intervals lying entirely inside $I_d$:
\[
 \mathcal D=\bigl\{0\le b<M:\ [b/M,(b+1)/M]\subset I_d\bigr\}.
\]
Since $M\ge4/\ell_d$, we have
\[
 |\mathcal D|\ge\ell_dM-2\ge\frac{\ell_dM}{2}.
\]
Thus the proportion of digits outside $\mathcal D$ is at most $1-\ell_d/2$. Write
\[
 \rho_d=1-\frac{\ell_d}{2}<1,
\]
so that $M-|\mathcal D|\le\rho_dM$. Let $N=\lceil\log(x+1)/\log M\rceil$, and write $m\le x$ as $m=\sum_{r=0}^{N-1}b_rM^r$, padding with leading zeros to obtain $N$ digits if necessary. Then
\[
 \Bigl\{\frac{m}{p^{L_d(r+1)}}\Bigr\}=\Bigl\{\frac m{M^{r+1}}\Bigr\}=\frac{b_r+\theta_r}{M},\qquad0\le\theta_r<1 .
\]
Whenever $b_r\in\mathcal D$, Lemma~\ref{lem2.2} ensures that the term of~\eqref{eq2.2} with index $L_d(r+1)$ equals $1$. Different digit positions contribute to different terms, so $v_p(C(m))\ge\#\{r:b_r\in\mathcal D\}$. If $v_p(C(m))\le A$, at most $A$ digits belong to $\mathcal D$. The number of such digit strings is at most
\[
 \sum_{r=0}^{\min(A,N)}\binom Nr|\mathcal D|^r(M-|\mathcal D|)^{N-r}
 \le\sum_{r=0}^{A}N^r\rho_d^{-r}(\rho_dM)^N
 \ll_{d,A}(1+N)^A(\rho_dM)^N.
\]
The factor $\rho_d^N$ in this bound decays geometrically with the number of digits $N$. Rewriting $(\rho_dM)^N$ as a power of $M^N$, we have
\[
 (\rho_dM)^N=(M^N)^{1+\log\rho_d/\log M}.
\]
We therefore set
\[
 \eta_p=-\frac{\log\rho_d}{\log M}
       =-\frac{\log\rho_d}{L_d\log p}>0.
\]
Now $M^N<M(x+1)$ and $N\ll_{d,p}1+\log x$ give~\eqref{eq3.1}. Since $L_d,\rho_d$ depend only on $d$, the exponent $\eta_p$ can be chosen uniformly for all $1\le j<d$.
\end{proof}

We retain the notation $L_d$ and $\rho_d$ above. The estimate at a fixed prime will be used in Section~\ref{sec6} to compare different greatest common divisor conditions.

We next estimate the coprime count restricted to integers with a small prime factor. For $m>1$, let $P^-(m)$ denote its least prime factor, with the convention $P^-(1)=+\infty$. We need the following short-interval bound.

\begin{lemma}\label{lem3.2}
There is an absolute constant $c_0$ such that, for every prime $p$ and real number $H$, the interval $[H,H+p)$ contains at most $c_0p/\log p$ integers with no prime factor less than $p$; that is, \(\#\{n\in[H,H+p):P^-(n)\ge p\}\ll \frac{p}{\log p}\).
\end{lemma}

\begin{proof}
This is a standard consequence of an upper-bound sieve~\cite[Theorem~2.2]{bib7}. We give a brief proof using Selberg's sieve. We may assume that $p$ is sufficiently large. Taking $z=p^{1/4}<p$, we obtain the bound
\[
 \ll \frac p{G(z)}+z^2,\qquad
 G(z)=\sum_{r\le z}\frac{\mu(r)^2}{\varphi(r)}\gg\log z,
\]
where $\mu$ and $\varphi$ are the M\"obius and Euler functions, respectively. The bound is therefore $O(p/\log p+p^{1/2})=O(p/\log p)$.
\end{proof}

We now use the carry interval \(I_d\) to restrict base-\(p\) digits. If the interval \([b/p,(b+1)/p]\) corresponding to a digit \(b\in\{0,1,\ldots,p-1\}\) lies entirely in \(I_d\), then the occurrence of this digit in the base-\(p\) expansion of \(m\) implies \(p\mid C(m)\), by Lemma~\ref{lem2.2}. Hence none of these digits can occur in an integer \(m\) when \(p\mid m\) and \(p\nmid C(m)\); we call them \emph{forbidden digits}. For sufficiently large \(p\), they constitute a fixed positive proportion of all \(p\) digits. Combining this restriction with the short-interval sieve bound in Lemma~\ref{lem3.2} allows us to estimate the coprime count for integers with small prime factors.

\begin{lemma}\label{lem3.3}
There exists $b_d>0$ such that, for $x\ge3$ and $0<\delta\le1/4$,
\begin{equation}
 \#\{m\le x:\gcd(m,C(m))=1,\ P^-(m)\le x^\delta\}
 \ll_d\frac x{\log x}\exp(-b_d/\delta).
 \label{eq3.2}
\end{equation}
\end{lemma}

\begin{proof}
Fix $4/\ell_d\le p\le x^\delta$ and count the integers $m\le x$ satisfying $P^-(m)=p$ and $p\nmid C(m)$. Put $D=\lfloor\log x/\log p\rfloor\ge4$ and write $m=\sum_{r=0}^{D}b_rp^r$, where $b_0=0$. If $[b_r/p,(b_r+1)/p]\subset I_d$, then $\{m/p^{r+1}\}\in I_d$, so $p\mid C(m)$. There are at least $\ell_dp-2\ge\ell_dp/2$ such digits. Consequently, when $p\nmid C(m)$, each digit has at most $\rho_dp$ possible values.

Apply this restriction to $b_2,\ldots,b_{D-1}$. The leading digit satisfies $b_D\le x/p^D$ and thus has at most $2x/p^D$ possible values. Once $b_2,\ldots,b_D$ are fixed, we have $m/p=H+b_1$, where $H=\sum_{r=2}^{D}b_rp^{r-1}$ and $0\le b_1<p$. Since $P^-(m)=p$, the integer $m/p$ has no prime factor less than $p$, so Lemma~\ref{lem3.2} bounds the number of choices for $b_1$ by $O(p/\log p)$. Combining these estimates gives $O_d\bigl(x\rho_d^{D-2}/(p\log p)\bigr)$. To express the geometric decay in exponential form, put
\[
 a_d=-\log\rho_d>0.
\]
Since $D-2\ge\log x/\log p-3$, we obtain
\begin{equation}
 \#\{m\le x:P^-(m)=p,\ p\nmid C(m)\}
 \ll\frac{x}{p^D}(\rho_dp)^{D-2}\frac{p}{\log p}
 \ll_d \frac x{p\log p}
   \exp\Bigl(-a_d\frac{\log x}{\log p}\Bigr).
 \label{eq3.3}
\end{equation}

To sum~\eqref{eq3.3}, first note that, for $2\le y\le x$,
\begin{equation}
 \sum_{p\le y}\frac{1}{p\log p}\exp\Bigl(-a_d\frac{\log x}{\log p}\Bigr)
 \ll_d\frac1{\log x}
 \exp\Bigl(-\frac{a_d\log x}{2\log y}\Bigr).
 \label{eq3.4}
\end{equation}
Indeed, put $g(t)=\exp(-a_d\log x/\log t)/(t\log t)$. Applying partial summation to $\sum_{p\le y}g(p)$ gives
\[\sum_{p\le y}g(p)
=\pi(y)g(y)-\int_2^y\pi(t)g'(t)\,dt.
\]
Use $\pi(t)\ll t/\log t$ and make the substitution
\[
 u=\frac{\log x}{\log t},\qquad
 \frac{dt}{t(\log t)^2}=-\frac{du}{\log x}.
\]
The absolute value of the integral is bounded by
\[
\begin{aligned}
&\ll_d \frac{1}{\log x}
\int_{\log x/\log y}^{\infty}
e^{-a_du}
\left(1+\frac {u}{\log x}+\frac{u^2}{\log x}\right)du\\
&\ll_d \frac{1}{\log x}
\int_{\log x/\log y}^{\infty}e^{-a_du/2}\,du\\
&\ll_d \frac{1}{\log x}
\exp\left(-\frac{a_d\log x}{2\log y}\right).
\end{aligned}
\]
On the other hand, since $\log x/\log y\ge1$
and $\log y\ge \log2$, the endpoint term satisfies
\[
\begin{aligned}
\pi(y)g(y)
\ll \frac{e^{-a_d\log x/\log y}}{(\log y)^2}\ll_d\frac{e^{-a_d\log x/2\log y}}{\log x}.
\end{aligned}
\]
Combining the integral and endpoint estimates proves~\eqref{eq3.4}.

Taking $y=x^\delta$, we see that the total contribution from $4/\ell_d\le p\le x^\delta$ is $O_d\bigl(xe^{-a_d/(2\delta)}/\log x\bigr)$.

For the finitely many primes $p<4/\ell_d$, Lemma~\ref{lem3.1} with $A=0$ gives a total contribution of $O_d(x^{1-\eta_0})$, where $\eta_0=a_d/(L_d\log 4/\ell_d)$. This range is nonempty only if $x^\delta\ge2$, in which case $1/\delta\le\log x/\log2$. Take
\[
 b_d=\min\{a_d/4,\ \eta_0\log2/4\},
\]
Then $e^{-b_d/\delta}\ge x^{-\eta_0/4}$, and hence $x^{1-\eta_0}\ll_d xe^{-b_d/\delta}/\log x$. Together with the estimate for the larger primes, this proves~\eqref{eq3.2}.
\end{proof}

\section{Simultaneous counting of carry conditions}\label{sec4}

When all prime factors are large, the number of carry conditions to be tested is bounded. We prove a simultaneous counting formula in which $R$ nontrivial conditions give the main-term weight $2^{-R}$, reflecting the measure $1/2$ of $Z$. We first approximate the indicator function of $Z$ by trigonometric polynomials, and then control the nonzero frequencies using exponential sum estimates over primes.

\subsection{Exponential sums and trigonometric approximation}\label{subsec4.1}

Write $e(t)=e^{2\pi it}$. We begin with the required exponential sum estimate.

\begin{lemma}[{\cite[Lemma~3.3]{bib4}}]\label{lem4.1}
Let $N\ge2$ be an integer and let $r_1\le r_2$ be positive integers. Put
\[
 f(u)=\alpha u+\sum_{r=r_1}^{r_2}\frac{\beta_r}{u^r},\qquad
 \alpha,\beta_{r_1},\ldots,\beta_{r_2}\in\mathbb R.
\]
Suppose that there exists $A\in[1,N^{1/6}]$ such that
\begin{equation}
 |\beta_{r_1}|\ge N^{r_1}A,\qquad |\beta_r/\beta_{r_1}|\le N^{(r-r_1)/2}\quad(r_1\le r\le r_2).
 \label{eq4.1}
\end{equation}
Then
\[
 \max_{I\subset(N,2N]}\Bigl|\sum_{p\in I}e(f(p))\Bigr|
 \ll_{r_2}N(\log N)^4\bigl(N^{-1/(3\cdot2^{\mathfrak j})}+A^{-1/10}\bigr),
 \quad
 \mathfrak j=3+\Bigl\lfloor\frac{\log\bigl(|\beta_{r_1}|/(AN^{r_1})\bigr)}{\log N}\Bigr\rfloor,
\]
where the maximum is over intervals $I$.
\end{lemma}

We also use Selberg's one-sided trigonometric approximation~(\cite[Lemma~4.1]{bib4}, \cite{bib13}). For every $K\in\mathbb N$ and every interval $I$ on the circle $\mathbb R/\mathbb Z$, there are trigonometric polynomials $S^\pm_{K,I}$ of degree at most $K$ such that
\[
 S^-_{K,I}\le\ind_I\le S^+_{K,I},\qquad
 \int_0^1S^\pm_{K,I}=|I|\pm\frac1{K+1},
\]
and the absolute value of the $n$th Fourier coefficient is at most $4/(|n|+1)$. Here $\ind_I$ denotes the indicator function of $I$.

Put $J_d=2d+1$. Decompose $Z$ into at most $J_d$ pairwise disjoint intervals and sum their upper and lower approximating polynomials. This gives polynomials $U,L$ of degree at most $K$ satisfying
\begin{equation}
 L\le\ind_Z\le U,\qquad U\ge0,\qquad
 \int_0^1U=\frac12+\frac{\theta_1J_d}{K+1},\quad
 \int_0^1(U-L)\le\frac{2J_d}{K+1},
 \label{eq4.2}
\end{equation}
where $0\le\theta_1\le1$, and all Fourier coefficients of $U,L$ are bounded in absolute value by $4J_d$. To test several conditions simultaneously, we use the following product inequality.

\begin{lemma}\label{lem4.2}
Let $L_\nu\le a_\nu\le U_\nu$ be real numbers for $1\le\nu\le R$, where $a_\nu\in\{0,1\}$ and $U_\nu\ge0$. Then
\[
 \prod_{\nu=1}^RU_\nu-\sum_{\nu=1}^R(U_\nu-L_\nu)\prod_{\mu\ne\nu}U_\mu
 \ \le\ \prod_{\nu=1}^Ra_\nu\ \le\ \prod_{\nu=1}^RU_\nu .
\]
\end{lemma}

\begin{proof}
The upper bound follows from $0\le a_\nu\le U_\nu$. A telescoping expansion gives
\[
 \prod_\nu U_\nu-\prod_\nu a_\nu
 =\sum_{\nu=1}^R(U_\nu-a_\nu)
   \prod_{\mu<\nu}a_\mu\prod_{\mu>\nu}U_\mu
 \le\sum_{\nu=1}^R(U_\nu-L_\nu)\prod_{\mu\ne\nu}U_\mu,
\]
Rearranging proves the lower bound.
\end{proof}

\subsection{The simultaneous counting formula}\label{subsec4.2}

Following~{\cite[Proposition~3]{bib4}}, we establish the simultaneous counting formula needed below. The factor $2^{-R}$ comes from the zero-frequency contribution of all $R$ conditions, while the nonzero frequencies are controlled by Lemma~\ref{lem4.1}.

\begin{lemma}\label{lem4.3}
Fix integers $k\ge2$, $H\ge1$ and real numbers $0<\varepsilon<1/8$, $0<\rho<\varepsilon/8$. Let $x$ be sufficiently large, and suppose that positive integers $M_1,\ldots,M_k$ and nonnegative integers $s_1,\ldots,s_k$ satisfy
\begin{equation}
 M_i\ge x^\varepsilon,\qquad
 \frac x{2^k}<\prod_{i=1}^kM_i\le2x,\qquad
 4dx^{1/(s_i+2)}\le M_i\le x^{(1-\rho)/(s_i+1)}.
 \label{eq4.3}
\end{equation}
Let $\mathcal Q$ be a set of prime tuples satisfying
\[
 q_i\in(M_i,2M_i],\qquad q_1<\cdots<q_k,\qquad
 x<q_1\cdots q_k\le2x
\]
together with additional restrictions such that, after any $k-1$ coordinates are fixed, the permitted range of the remaining coordinate is a union of at most $t$ intervals. Put $R=s_1+\cdots+s_k$, and suppose that $H\ge4J_dR$. Then there exists $\xi>0$, depending only on $k,\varepsilon,\rho$, such that
\begin{equation}
 \#\bigl\{\bm q\in\mathcal Q:\gcd(m,C(m))=1\bigr\}
 =2^{-R}\Bigl(1+O_d\Bigl(\frac RH\Bigr)\Bigr)|\mathcal Q|
  +O_{d,k,\varepsilon,\rho,H}(tx^{1-\xi}),
 \label{eq4.4}
\end{equation}
where $m=q_1\cdots q_k$.
\end{lemma}

\begin{proof}
\noindent\textit{Truncating the carry conditions.}
Put $m=q_1\cdots q_k$. Since the $q_i$ are distinct, coprimality is equivalent to every $v_{q_i}(C(m))$ being $0$. In~\eqref{eq2.2}, with $p=q_i$, the term with $r=1$ is $F(0)=0$. For $r\ge s_i+2$, condition~\eqref{eq4.3} gives
\[
 \frac{m}{q_i^r}\le\frac{2x}{M_i^{s_i+2}}\le\frac{2}{(4d)^{s_i+2}}<\frac1d,
\]
so the corresponding terms also equal $0$. Since all remaining terms are nonnegative, we have precisely
\[
 \gcd(m,C(m))=1\iff
 \Bigl\{\frac m{q_i^{s+1}}\Bigr\}\in Z\qquad(1\le i\le k,\ 1\le s\le s_i).
\]
There are $R$ conditions. The inequalities $M_i\ge x^\varepsilon$ and $M_i^{s_i+1}\le x$ give $R\le k/\varepsilon$. The conclusion is immediate when $R=0$, so we assume that $R\ge1$.

\medskip\noindent\textit{The zero-frequency main term.}
Take $K=H$ in~\eqref{eq4.2} and apply Lemma~\ref{lem4.2} to the $R$ conditions. Lower and upper bounds for the desired count are, respectively,
\[
 \sum_{\bm q\in\mathcal Q}\Bigl(\prod_{\nu}U(t_\nu)-\sum_\nu(U-L)(t_\nu)\prod_{\mu\ne\nu}U(t_\mu)\Bigr)
 \quad\text{and}\quad
 \sum_{\bm q\in\mathcal Q}\prod_{\nu}U(t_\nu)
\]
where the $t_\nu$ run through $m/q_i^{s+1}$ for $1\le i\le k$ and $1\le s\le s_i$. Expanding the trigonometric polynomials, we obtain zero-frequency contributions equal to $|\mathcal Q|$ times
\[
 \Bigl(\int_0^1U\Bigr)^R-R\int_0^1(U-L)\Bigl(\int_0^1U\Bigr)^{R-1}
 \quad\text{and}\quad
 \Bigl(\int_0^1U\Bigr)^R.
\]
By~\eqref{eq4.2} and $H\ge4J_dR$, we have $(1+2J_d/(H+1))^R\le e^{2J_dR/H}=1+O_d(R/H)$. Thus both zero-frequency coefficients are $2^{-R}(1+O_d(R/H))$.

\medskip\noindent\textit{The nonzero-frequency error.}
The nonzero-frequency terms are finitely many exponential sums with bounded coefficients, and their number depends only on $d,H,R$. Each sum has the form
\[
 \sum_{\bm q\in\mathcal Q}e(\Phi(\bm q)),\qquad
 \Phi(\bm q)=\sum_{i=1}^k\sum_{s=1}^{s_i}\lambda_{i,s}\frac{m}{q_i^{s+1}},
\]
where the $\lambda_{i,s}$ are integers, not all zero, with $|\lambda_{i,s}|\le H$. Choose $h$ such that the $h$th row is nonzero, and put $r_1=\min\{s:\lambda_{h,s}\ne0\}$. Fix the other primes and write $P=\prod_{i\ne h}q_i$. The terms with $i\ne h$ are linear in $q_h$, so there is a real number $\alpha$, independent of $q_h$, such that
\[
 \Phi(\bm q)=\alpha q_h+\sum_{s=r_1}^{s_h}\frac{P\lambda_{h,s}}{q_h^{s}}
\]
To apply Lemma~\ref{lem4.1}, take $N=M_h$, $\beta_s=P\lambda_{h,s}$, and $A=x^a$, where $a=\min\{\rho/2,\varepsilon/12\}$. For a nonempty permitted section, the inequalities $m>x$, $q_h\le2M_h$, and $M_h^{s_h+1}\le x^{1-\rho}$ give
\[
 P>\frac{x}{2M_h}\ge\frac12x^\rho M_h^{s_h},
\]
Since $r_1\le s_h$, for sufficiently large $x$ we have $|\beta_{r_1}|\ge\tfrac12x^\rho M_h^{r_1}\ge M_h^{r_1}A$. Also, $A\le x^{\varepsilon/12}\le M_h^{1/6}$, and, for $s>r_1$, we have $|\beta_s/\beta_{r_1}|\le H\le M_h^{(s-r_1)/2}$. This verifies~\eqref{eq4.1}. Moreover, $|\beta_{r_1}|\le HP\le2Hx/M_h$ and $AN^{r_1}\ge M_h$ imply that, for sufficiently large $x$,
\[
 \mathfrak j\le3+\frac{\log(2Hx)}{\log M_h}-2\le 3+\frac1\varepsilon=:\mathfrak j_0
\]
Therefore, on each permitted interval,
\[
 \sum_{q_h}e(\Phi(\bm q))\ll_{k,\varepsilon}M_h(\log x)^4\bigl(x^{-\varepsilon/(3\cdot2^{\mathfrak j_0})}+x^{-a/10}\bigr)\ll M_hx^{-\xi},
\]
where we may take $\xi=\tfrac12\min\{\varepsilon/(3\cdot2^{\mathfrak j_0}),a/10\}$.
Summing first over the at most $t$ permitted intervals
and then over the other $k-1$ prime coordinates,
and using $\prod_{i=1}^k M_i\le 2x$, we obtain
\[
\sum_{\boldsymbol q\in\mathcal Q}
e(\Phi(\boldsymbol q))
\ll
tM_hx^{-\xi}\prod_{i\ne h}M_i
\ll tx^{1-\xi}.
\]
Since $R\le k/\varepsilon$ and $H$ is fixed,
the number of nonzero frequencies depends only on $k,\varepsilon,H$.
Summing over these frequencies proves~\eqref{eq4.4}.
\end{proof}

\section{An asymptotic formula for the coprime count}\label{sec5}

\subsection{The main-term constant and its convergence}\label{subsec5.1}

We now use Lemma~\ref{lem4.3} to compute the main term for each number of distinct prime factors. Away from the critical regions near \(x^{1/r}\), a prime factor \(q=x^u\) contributes \(\lfloor1/u\rfloor-1\) nontrivial carry conditions; see~\eqref{eq5.7}. The main weight in Lemma~\ref{lem4.3} can therefore be written as

\[
2^{-R}=\prod_{i=1}^{r}2^{\,1-\lfloor1/u_i\rfloor}.
\]

On the other hand, the substitution \(y=x^u\) gives

\[
\frac{dy}{\log y}=\frac{x^u}{u}\,du.
\]

This leads us to introduce the weight function

\begin{equation}
h(u)=\frac{2^{\,1-\lfloor1/u\rfloor}}{u}
\qquad(0<u\le1),
\label{eq5.1}
\end{equation}

and put \(\mathfrak c_1=1\). For \(r\ge2\), define

\begin{equation}
\mathfrak c_r=
\frac1{r!}
\int_{\substack{u_1,\ldots,u_{r-1}>0\\
u_1+\cdots+u_{r-1}<1}}
\prod_{i=1}^{r-1}h(u_i)\,
h(1-u_1-\cdots-u_{r-1})
\,du_1\cdots du_{r-1},
\label{eq5.2}
\end{equation}

and set

\[
\mathfrak c=\sum_{r\ge1}\mathfrak c_r.
\]

This is the same constant as that of Ford and Konyagin. Indeed, if \(1/(n+1)<u\le1/n\), then

\[
0<h(u)=\frac{2^{1-n}}u
\le 2^{1-n}(n+1)\le2.
\]

Since the region of integration has volume \(1/(r-1)!\), we obtain

\[
0<\mathfrak c_r
\le\frac{2^r}{r!(r-1)!}.
\]

Hence the series defining \(\mathfrak c\) converges absolutely.

\subsection{Counting with a fixed number of prime factors}\label{subsec5.2}

Let $\omega(m)$ denote the number of distinct prime factors of $m$, and put
\[
 N_k(x)=\#\{x<m\le2x:\omega(m)=k,\ \gcd(m,C(m))=1\}.
\]
We show that $\cc_k$ is the main-term coefficient for fixed $k$.

\begin{lemma}\label{lem5.1}
For each fixed integer $k\ge1$,
\begin{equation}
 N_k(x)\sim\cc_k\frac x{\log x}\qquad(x\to\infty).
 \label{eq5.3}
\end{equation}
\end{lemma}

\begin{proof}
When $k=1$, the integers being counted are prime powers. For a prime $q>d$, Lucas's theorem gives $C(q)\equiv\binom dj\not\equiv0\pmod q$, so every such $q$ satisfies the coprimality condition. There are only $O(\sqrt x)$ prime powers with exponent at least $2$ in $(x,2x]$. The prime number theorem therefore gives $N_1(x)\sim x/\log x=\cc_1x/\log x$.

Fix \(k\ge2\). The proof has two steps. We first exclude integers with small prime factors, square factors, or prime factors in critical regions near \(x^{1/r}\), apply Lemma~\ref{lem4.3} to the remaining prime tuples, and evaluate their main term using the prime number theorem. We then estimate the excluded sets separately and show that they do not affect the final main-term coefficient. To apply Lemma~\ref{lem4.3}, we place the prime factors in dyadic intervals away from the critical points \(x^{1/r}\), so that the number of carry conditions is constant on each block and the required separation conditions hold.

\medskip\noindent\textit{Truncation and critical windows.}
Take
\[
0<\varepsilon<\frac18,\qquad
\rho=\varepsilon^4,\qquad
V=\left\lceil\frac2\varepsilon\right\rceil+2.
\]
Define the critical windows by
\begin{equation}
W(x)=\bigcup_{r=1}^{V}\mathcal W_r(x),
\qquad
\mathcal W_r(x)=
\left(x^{(1-\rho)/r},\,8dx^{1/r}\right].
\label{eq5.4}
\end{equation}
The choice of \(V\) covers all critical exponents that can occur in the subsequent decomposition, while the choice of \(\rho\) controls the window error as \(\varepsilon\to0\).

Let \(G(x)=G_{k,\varepsilon}(x)\) count the integers satisfying
\[
x<m\le2x,\qquad
m=q_1\cdots q_k,\qquad
q_1<\cdots<q_k,
\]
\[
q_i>x^\varepsilon,\qquad
q_i\notin W(x),\qquad
\gcd(m,C(m))=1
\]
where the \(q_i\) are primes. Define the corresponding truncated integral constant by
\begin{equation}
\cc_{k,\varepsilon}
=
\frac1{k!}
\int_{\substack{v_i\ge\varepsilon\\
v_1+\cdots+v_k=1}}
\prod_{i=1}^kh(v_i)\,
dv_1\cdots dv_{k-1},
\label{eq5.5}
\end{equation}
where \(v_k=1-v_1-\cdots-v_{k-1}\).

The omitted part of the simplex has volume \(O_k(\varepsilon)\), and \(0<h(u)\le2\). Hence
\[
0\le\mathfrak c_k-\mathfrak c_{k,\varepsilon}
\ll_k\varepsilon.
\]

\medskip\noindent\textit{Dyadic decomposition and simultaneous detection.}
For each prime factor, choose a power of $2$, denoted by $M_i$, such that $M_i<q_i\le2M_i$. For a given $\bm M=(M_1,\ldots,M_k)$, let $\mathcal T(\bm M)$ be the set of real vectors satisfying
\[
 y_i\in(M_i,2M_i],\quad y_1<\dots<y_k,\quad x<y_1\cdots y_k\le2x,\quad y_i>x^\varepsilon,\quad y_i\notin W(x)
\]
Let $\mathcal Q(\bm M)$ be the set of prime tuples in $\mathcal T(\bm M)$. After the other coordinates are fixed, the interval, ordering, product, and lower-bound restrictions define an interval; the window restriction removes at most $V$ intervals from it. Thus we may take $t=V+1$ in Lemma~\ref{lem4.3}. If $\mathcal T(\bm M)$ is nonempty, then
\[
 M_i\ge x^\varepsilon/2\ge x^{\varepsilon/2},\qquad x/2^k<\prod_iM_i\le2x .
\]
Since $k\ge2$ and all other coordinates exceed $x^\varepsilon$, we have $M_i<2x^{1-\varepsilon}<x$. Put $s_i=\lfloor\log x/\log M_i\rfloor-1$. Then $s_i\ge0$, $x^{1/(s_i+2)}<M_i\le x^{1/(s_i+1)}$, and $s_i+2\le2/\varepsilon+1\le V$. The window restriction further ensures that
\begin{equation}
 4dx^{1/(s_i+2)}\le M_i\le x^{(1-\rho)/(s_i+1)}.
 \label{eq5.6}
\end{equation}
If the left inequality failed, the whole interval $(M_i,2M_i]$ would lie in $\mathcal W_{s_i+2}$; if the right inequality failed, it would lie in $\mathcal W_{s_i+1}$. Both cases contradict the nonemptiness of $\mathcal T(\bm M)$. Since $\rho=\varepsilon^4<\varepsilon/16$, Lemma~\ref{lem4.3} applies with parameters $(k,\varepsilon/2,\rho)$. Moreover, for $\bm y\in\mathcal T(\bm M)$, we have $y_i>M_i>x^{1/(s_i+2)}$ and $y_i\le2x^{(1-\rho)/(s_i+1)}<x^{1/(s_i+1)}$, so
\begin{equation}
 s_i=\Bigl\lfloor\frac{\log x}{\log y_i}\Bigr\rfloor-1 .
 \label{eq5.7}
\end{equation}

Take $H=\lceil16J_dk\varepsilon^{-4}\rceil$. Since $R=\sum s_i\le2k/\varepsilon$, we have $H\ge4J_dR$ and $R/H\ll_d\varepsilon^3$. There are at most $O_k((\log x)^k)$ nonempty blocks, so Lemma~\ref{lem4.3} gives
\begin{equation}
 G(x)=\bigl(1+O_d(\varepsilon^3)\bigr)\sum_{\bm M}2^{-R(\bm M)}|\mathcal Q(\bm M)|+o(x/\log x).
 \label{eq5.8}
\end{equation}

\medskip\noindent\textit{Summation over primes.}
Each coordinate section is a union of at most $V+1$ intervals,
and $M_i\ge x^{\varepsilon/2}$.
The prime number theorem with its classical error term gives, uniformly for every interval
$I\subset(M_i,2M_i]$,
\[
\#\{p\in I:p\text{ is prime}\}
=
\int_I\frac{dt}{\log t}
+
O\left(M_i e^{-c_0\sqrt{\log M_i}}\right),
\]
where $c_0>0$ is an absolute constant.
Fixing the remaining coordinates and applying this estimate successively to each coordinate,
we use $\prod_iM_i\le 2x$ to obtain
\[
|\mathcal Q(\bm M)|
=
\int_{\mathcal T(\bm M)}
\frac{dy_1\cdots dy_k}
{\log y_1\cdots\log y_k}
+
O_{k,\varepsilon}
\left(xe^{-c'\sqrt{\varepsilon\log x}}\right),
\]
for a suitable constant $c'>0$.

Write $\mathcal T=\bigcup_{\bm M}\mathcal T(\bm M)$. By~\eqref{eq5.7},
on each nonempty block $\mathcal T(\bm M)$ we have
\[
2^{-R(\bm M)}
=
\prod_{i=1}^k
2^{1-\lfloor\log x/\log y_i\rfloor}.
\]
There are at most
$O_k((\log x)^k)$ nonempty blocks, so the total error from the prime number theorem is still
\[
O_{k,\varepsilon}\left(
x(\log x)^k e^{-c'\sqrt{\varepsilon\log x}}
\right)
=o(x/\log x).
\]
The blocks have disjoint interiors. Combining the integrals
and using~\eqref{eq5.8}, we obtain
\begin{equation}
G(x)
=
\bigl(1+O_d(\varepsilon^3)\bigr)
\int_{\mathcal T}
\prod_{i=1}^k
\frac{2^{1-\lfloor\log x/\log y_i\rfloor}}
{\log y_i}\,d\boldsymbol y
+
o(x/\log x).
\label{eq5.9}
\end{equation}

\medskip\noindent\textit{Evaluation of the main-term integral.}
Make the substitution $y_i=x^{u_i}$. By the definition of $h$,
\[
\frac{2^{1-\lfloor\log x/\log y_i\rfloor}}
{\log y_i}\,dy_i
=x^{u_i}h(u_i)\,du_i.
\]
Remove the ordering condition and multiply by $1/k!$. Next put
$u_i=tv_i$, where $\sum_{i=1}^k v_i=1$.
The product restriction gives
$1<t\le 1+\log 2/\log x$,
and the absolute value of the Jacobian is $t^{k-1}$.
Temporarily omitting the window restrictions, we obtain the integral
\[
\frac1{k!}
\int_1^{1+\log 2/\log x}
x^t t^{k-1}
\int_{\substack{v_i>\varepsilon/t\\
                \sum_{i=1}^k v_i=1}}
\prod_{i=1}^k h(tv_i)
\,dv_1\cdots dv_{k-1}\,dt.
\]

We compare the inner integral with the definition of $\cc_{k,\varepsilon}$.
Since $\varepsilon$ is fixed, $h$ is bounded on
$[\varepsilon/2,1]$, has only finitely many jump discontinuities,
and is piecewise $C^1$ elsewhere.
If the interval between $v_i$ and $tv_i$ crosses a jump point
$\alpha$, then
\[
\frac{\alpha}{t}\le v_i\le\alpha.
\]
This interval has length $O(t-1)$.
Since there are only finitely many jump points, the corresponding regions in the simplex
have total volume $O_{k,\varepsilon}(t-1)$.
Similarly, replacing the boundary condition $v_i>\varepsilon/t$
by $v_i\ge\varepsilon$ changes the region by a thin layer
of the same order of volume.

On the remaining region, $v_i$ and $tv_i$
lie in the same interval of smoothness of $h$.
The mean value theorem and the boundedness of $h$ give
\[
\prod_{i=1}^k h(tv_i)
=
\prod_{i=1}^k h(v_i)
+O_{k,\varepsilon}(t-1).
\]
Using also $t^{k-1}=1+O_k(t-1)$,
we obtain
\[
t^{k-1}
\int_{\substack{v_i>\varepsilon/t\\
                \sum_{i=1}^k v_i=1}}
\prod_{i=1}^k h(tv_i)\,d\boldsymbol v
=
k!\cc_{k,\varepsilon}
+O_{k,\varepsilon}(t-1),
\]
where $d\boldsymbol v=dv_1\cdots dv_{k-1}$.
Thus, without the window restrictions, the integral equals
\[
\begin{aligned}
&\frac1{k!}
\int_1^{1+\log 2/\log x}
x^t\bigl(k!\cc_{k,\varepsilon}
+O_{k,\varepsilon}(t-1)\bigr)\,dt\\
&\qquad=
\cc_{k,\varepsilon}\frac{x}{\log x}
+O_{k,\varepsilon}\left(\frac{x}{\log^2 x}\right).
\end{aligned}
\]

It remains to estimate the effect of the window restrictions.
By definition, the length of $W_r$ in the $u$-variable is
\[
\frac{\rho}{r}+\frac{\log(8d)}{\log x}.
\]
For fixed $t$, the region of the simplex in which a given coordinate
$tv_i$ belongs to this window has volume
at most $O_k(1)$ times the window length.
Summing over all coordinates and windows, we find that the excluded region
has volume at most
\[
\begin{aligned}
&\ll_k\sum_{r\le V}
\left(\frac{\rho}{r}
+\frac{\log(8d)}{\log x}\right)\\
&\ll_d \rho\log(2/\varepsilon)
+\frac{1}{\varepsilon\log x},
\end{aligned}
\]
where we used $V\asymp1/\varepsilon$.
Together with $h\le2$ and the above range of integration in $t$,
this bounds the change in the integral by
\[
\begin{aligned}
&\ll_{d,k}\frac{x}{\log x}
\left(
\rho\log(2/\varepsilon)
+\frac{1}{\varepsilon\log x}
\right)\\
&=
O_{d,k}\left(
\varepsilon^3\frac{x}{\log x}
\right)+o(x/\log x),
\end{aligned}
\]
Here we used $\rho=\varepsilon^4$
and let $x\to\infty$ with $k,\varepsilon$ fixed.

Combining these estimates with the preceding prime-summation formula,
and noting that $\cc_{k,\varepsilon}\ll_k1$, we obtain
\begin{equation}
G(x)=\frac{x}{\log x}
\left(\cc_{k,\varepsilon}
+O_{d,k}(\varepsilon^3)+o(1)\right).
\label{eq5.10}
\end{equation}

\medskip\noindent\textit{Exceptional sets and the limiting argument.}
By definition, $G(x)\le N_k(x)$. We now consider only integers counted by $N_k(x)$ but not by $G(x)$. They all satisfy $\gcd(m,C(m))=1$ and belong to at least one of the following three classes.

(a) Integers with a prime factor $\le x^\varepsilon$. Applying Lemma~\ref{lem3.3} with upper bound $2x$ and $\delta=2\varepsilon$, we obtain at most $O_d\bigl(e^{-b_d/(2\varepsilon)}x/\log x\bigr)$ such integers.

(b) Integers divisible by $p^2$ for some $p>x^\varepsilon$. There are at most $\sum_{p>x^\varepsilon}2x/p^2\ll x^{1-\varepsilon}$ such integers.

(c) Squarefree integers all of whose prime factors exceed $x^\varepsilon$, but with a prime factor $q\in \mathcal W_r$. Since the other prime factors exceed $x^\varepsilon$, we have $q\le2x^{1-\varepsilon}<x^{1-\rho}$, so only $r\ge2$ need be considered. Fix $q\in \mathcal W_r$ and another $k-2$ prime factors $q_2,\ldots,q_{k-1}$, and put $a=qq_2\cdots q_{k-1}$. The last prime factor satisfies $x^\varepsilon<q_k\le2x/a$, so there are at most $O\bigl(x/(\varepsilon a\log x)\bigr)$ choices for it. By Mertens' theorem,
\[
 \sum_{x^\varepsilon<p\le2x}\frac1p\ll\log(2/\varepsilon),\qquad
 \sum_{q\in \mathcal W_r}\frac1q=-\log(1-\rho)+O_{d,\varepsilon}\Bigl(\frac1{\log x}\Bigr)\ll\rho+O_{d,\varepsilon}\Bigl(\frac1{\log x}\Bigr)
\]
where the second estimate holds for $r\le V$. Thus the number of integers in this class is
\[
 \ll_{k}\frac x{\varepsilon\log x}
 (\log(2/\varepsilon))^{k-2}
 \sum_{r=2}^{V}\sum_{q\in \mathcal W_r}\frac1q
 \ll_{d,k}\varepsilon^2(\log(2/\varepsilon))^{k-2}\frac x{\log x}
 +o\Bigl(\frac x{\log x}\Bigr)
\]

Combining the three estimates with~\eqref{eq5.10} and using $|\cc_k-\cc_{k,\varepsilon}|\ll_k\varepsilon$, we obtain
\begin{equation}
 \limsup_{x\to\infty}\Bigl|\frac{N_k(x)\log x}{x}-\cc_k\Bigr|
 \ll_{d,k}\varepsilon+e^{-b_d/(2\varepsilon)}
 +\varepsilon^2(\log(2/\varepsilon))^{k-2}.
 \label{eq5.11}
\end{equation}
Finally, letting $\varepsilon\to0^+$ proves~\eqref{eq5.3}.
\end{proof}

\subsection{Summing over the number of prime factors}\label{subsec5.3}

\begin{proof}[Proof of Theorem~\ref{thm1.2} when $e=1$]
To sum over $k$, we must control the tail for large $k$. Fix an integer $K\ge8$. If $x<m\le2x$ and $\omega(m)>K$, then $P^-(m)\le(2x)^{1/(K+1)}$. Applying Lemma~\ref{lem3.3} with upper bound $2x$ and $\delta=2/(K+1)\le1/4$, we obtain
\[
 \#\{x<m\le2x:\gcd(m,C(m))=1,\ \omega(m)>K\}
 \ll_d e^{-b_dK/2}\frac x{\log x}.
\]
Apply Lemma~\ref{lem5.1} for $1\le k\le K$ and use the convergence of the series defining the constant to get
\[
 \limsup_{x\to\infty}
 \Bigl|\frac{B_{d,j}(2x)-B_{d,j}(x)}{x/\log x}-\cc\Bigr|
 \ll_d e^{-b_dK/2}+\sum_{k>K}\cc_k .
\]
Letting $K\to\infty$ gives $B_{d,j}(2x)-B_{d,j}(x)\sim\cc x/\log x$. Decomposing $(\sqrt y,y]$ into dyadic intervals and summing, while bounding the part up to $\sqrt y$ by $B_{d,j}(\sqrt y)\le\sqrt y$, yields $B_{d,j}(y)\sim\cc y/\log y$.
\end{proof}

\section{Completion of the main theorems and applications}\label{sec6}

With the coprime asymptotic established, we complete the proofs of the main theorems by using the fixed-prime valuation estimates from Section~\ref{sec3} to compare different greatest common divisor conditions.

\subsection{Recovering the original count: proof of Theorem~\ref{thm1.1}}\label{subsec6.1}

We first compare $T_{d,j}(x)$ from Lemma~\ref{lem2.1} with the coprime count $B_{d,j}(x)$.

\begin{lemma}\label{lem6.1}
There exists $\eta_d>0$ such that, for every $1\le j<d$,
\begin{equation}
 T_{d,j}(x)=B_{d,j}(x)+O_d(x^{1-\eta_d}).
 \label{eq6.1}
\end{equation}
\end{lemma}

\begin{proof}
If $p\nmid d$, both counting conditions at $p$ are equivalent to $\min(v_p(m),v_p(C(m)))=0$. They can differ only at primes $p\mid d$. For such a prime, write $a=v_p(d)$, $b=v_p(m)$, and $v=v_p(C(m))$. The two conditions are $\min(a+b,v)=a$ and $\min(b,v)=0$, respectively. When $v>a$, both are equivalent to $b=0$. Thus the symmetric difference of the two sets being counted is contained in
\[
 \bigcup_{p\mid d}\{m\le x:v_p(C(m))\le v_p(d)\},
\]
By Lemma~\ref{lem3.1}, this union has at most $O_d(x^{1-\eta_d})$ elements, where we may take $\eta_d=\tfrac12\min_{p\mid d}\eta_p$.
\end{proof}

\begin{proof}[Proof of Theorem~\ref{thm1.1}]
When $d=1$, the endpoints $k=0,n$ always satisfy the condition. If some $0<k<n$ also satisfied $\gcd\bigl(n,\binom nk\bigr)=1$, then $n\mid k\binom nk$ would imply $n\mid k$, a contradiction. Hence $A_1(n)=2$, and therefore $S_1(x)=2\lfloor x\rfloor$.

For $d\ge2$, applying Lemma~\ref{lem2.1}, Lemma~\ref{lem6.1}, and the coprime asymptotic~\eqref{eq1.4} in turn gives
\[
 S_d(x)=\sum_{j=1}^{d-1}B_{d,j}(x/d)+O_d\bigl(x^{1-\eta_d}\bigr)
 \sim (d-1)\,\cc\,\frac{x/d}{\log(x/d)}\sim\frac{d-1}{d}\,\cc\,\frac{x}{\log x}.
 \qedhere
\]
\end{proof}

\subsection{Counting a general fixed greatest common divisor}\label{subsec6.2}

The next proposition reduces the count for an arbitrary fixed greatest common divisor to a coprime count with different parameters.

\begin{proposition}\label{prop6.2}
For fixed integers $d\ge2$, $1\le j<d$, and $e\ge1$, there exists $\eta>0$, depending only on $d,e$, such that
\begin{equation}
 B_{d,j}(x;e)=B_{de,je}(x/e)+O_{d,e}(x^{1-\eta}).
 \label{eq6.2}
\end{equation}
\end{proposition}

\begin{proof}
When $e=1$, the assertion is an identity and we may take $\eta=1/2$. Suppose henceforth that $e>1$. A greatest common divisor equal to $e$ implies $e\mid m$, so write $m=et$ and put
\[
 C'(t)=\binom{(de)t}{(je)t}.
\]
This substitution changes the range to $t\le x/e$. It therefore suffices to compare the conditions $\gcd(et,C'(t))=e$ and $\gcd(t,C'(t))=1$.

Comparing prime by prime, a difference can occur only when $p\mid e$. If $p\nmid e$, both conditions require $\min(v_p(t),v_p(C'(t)))=0$. If $p\mid e$, write $a=v_p(e)$, $b=v_p(t)$, and $v=v_p(C'(t))$. The two conditions are $\min(a+b,v)=a$ and $\min(b,v)=0$, respectively. When $v>a$, both are equivalent to $b=0$. Hence the symmetric difference of the corresponding sets is contained in
\[
 \bigcup_{p\mid e}\{t\le x/e:v_p(C'(t))\le v_p(e)\}.
\]
Apply Lemma~\ref{lem3.1} with parameters $(de,je)$. Its proof gives the exponent
\[
 \eta_p=-\frac{\log\rho_{de}}{L_{de}\log p}
       =\frac{a_{de}}{L_{de}\log p}>0.
\]
Since there are only finitely many primes $p\mid e$, we may take $\eta=\tfrac12\min_{p\mid e}\eta_p>0$. The size of the symmetric difference is at most
\[
 \ll_{d,e}\sum_{p\mid e}x^{1-\eta_p}(1+\log x)^{v_p(e)}
 \ll_{d,e}x^{1-\eta}.
\]
Each fixed power of the logarithm can be absorbed into $x^{\eta_p/2}$. The constants in Lemma~\ref{lem3.1} depend only on $de,p,v_p(e)$, so the resulting estimate is uniform for $1\le j<d$. This proves~\eqref{eq6.2}.

\end{proof}

\begin{proof}[Proof of Theorem~\ref{thm1.2} for general $e$]
The coprime case~\eqref{eq1.4} was proved in the preceding section. Applying it with parameters $(de,je)$ gives
\[
 B_{de,je}(x/e)\sim\cc\frac{x/e}{\log(x/e)}
 \sim\frac{\cc}{e}\frac x{\log x}.
\]
Proposition~\ref{prop6.2} and $x^{1-\eta}=o(x/\log x)$ now imply~\eqref{eq1.3}.
\end{proof}

\subsection{Applications to primes, composites, and fixed ratios}\label{subsec6.3}

We may also split the count according to whether $n/d$ is prime. Put
\[
 \mathscr P_d(x)=\sum_{p\le x}A_d(dp),\qquad
 \mathscr C_d(x)=\sum_{\substack{2\le m\le x\\m\text{ composite}}}A_d(dm),
\]
where $p$ runs over primes.

\begin{corollary}\label{cor6.3}
For each fixed integer $d\ge2$,
\begin{equation}
 \mathscr P_d(x)\sim(d-1)\frac x{\log x},\qquad
 \mathscr C_d(x)\sim(d-1)(\cc-1)\frac x{\log x}.
 \label{eq6.3}
\end{equation}
\end{corollary}

\begin{proof}[Proof of Corollary~\ref{cor6.3}]
By the proof of Lemma~\ref{lem2.1}, $A_d(dp)$ is the number of $1\le j<d$ satisfying $\gcd\bigl(dp,\binom{dp}{jp}\bigr)=d$. If $p>d$, Lucas's theorem gives $\binom{dp}{jp}\equiv\binom dj\not\equiv0\pmod p$, so $\gcd\bigl(p,\binom{dp}{jp}\bigr)=1$.

The proof of Lemma~\ref{lem6.1} shows that, for all $1\le j<d$, this coprimality condition and the condition that the greatest common divisor equals $d$ can differ only on a set of $m\le x$ of size $O_d(x^{1-\eta_d})$. Thus $A_d(dp)=d-1$ for all but at most $O_d(x^{1-\eta_d})$ primes. It follows that
\[
 \mathscr P_d(x)=(d-1)\pi(x)+O_d(x^{1-\eta_d})
 \sim(d-1)\frac x{\log x}.
\]

Also, since $A_d(n)=0$ whenever $d\nmid n$, we have
\[
 S_d(dx)=A_d(d)+\mathscr P_d(x)+\mathscr C_d(x).
\]
By Theorem~\ref{thm1.1}, the left-hand side is asymptotic to $(d-1)\cc x/\log x$. Subtracting the prime contribution and the fixed term $A_d(d)$ gives the asymptotic formula for the composite part.
\end{proof}

\begin{remark}\label{rem6.4}
For $d=2$, there is an exact relation between the two counts. Write $B(Y)=\#\{m\le Y:\gcd\bigl(m,\binom{2m}{m}\bigr)=1\}$. Then
\begin{equation}
 S_2(x)=B(x/2)+\Bigl\lfloor\log_2\frac x2\Bigr\rfloor\qquad(x\ge2).
 \label{eq6.4}
\end{equation}
Indeed, only $j=1$ occurs, so $S_2(x)=T_{2,1}(x/2)$. By Kummer's theorem, $v_2\bigl(\binom{2m}{m}\bigr)=s_2(m)$, where $s_2(m)$ is the binary digit sum of $m$. If $m$ is odd, then $s_2(m)\ge1$, and $\gcd\bigl(2m,\binom{2m}{m}\bigr)=2$ is equivalent to $\gcd\bigl(m,\binom{2m}{m}\bigr)=1$. If $m$ is even, the latter condition fails, whereas the former holds if and only if $s_2(m)=1$, that is, if and only if $m$ is a positive integral power of $2$. There are $\lfloor\log_2(x/2)\rfloor$ such additional values of $m$, proving~\eqref{eq6.4}. Thus Theorem~\ref{thm1.1} for $d=2$ is equivalent to~\eqref{eq1.1}.
\end{remark}

One can also fix the ratio $j/d$ of the lower to the upper index and consider the greatest common divisor of the binomial coefficient and the upper index itself. The count has different orders of magnitude according to whether $d\mid e$ or $d\nmid e$.

\begin{corollary}\label{cor6.5}
Fix integers $d\ge2$, $1\le j<d$, and $e\ge1$. If $d\mid e$, then
\[
 \#\Bigl\{n\le x:d\mid n,\ \gcd\Bigl(n,\binom n{jn/d}\Bigr)=e\Bigr\}
 \sim\frac{\cc}{e}\frac{x}{\log x}\qquad(x\to\infty).
\]
If $d\nmid e$, there exists $\eta>0$, depending only on $d,e$, such that the above count is $O_{d,e}(x^{1-\eta})$.
\end{corollary}

\begin{proof}
First suppose that $d\mid e$. Since the greatest common divisor is $e$, we have $e\mid n$ and may write $n=et$. Then $d\mid n$ holds automatically. Moreover, $e\ge2$, $je/d$ is an integer, and $1\le je/d<e$, so the desired count is exactly $T_{e,je/d}(x/e)$. Lemma~\ref{lem6.1} and the coprime asymptotic~\eqref{eq1.4} give
\[
 T_{e,je/d}(x/e)=B_{e,je/d}(x/e)+O_e(x^{1-\eta_e})
 \sim\frac{\cc}{e}\frac{x}{\log x}.
\]

If $d\nmid e$, choose a prime $p\mid d$ such that $v_p(e)<v_p(d)$. Write $n=dm$ and retain the notation $C(m)=\binom{dm}{jm}$. If $\gcd(dm,C(m))=e$, then $v_p(dm)\ge v_p(d)>v_p(e)$ forces $v_p(C(m))=v_p(e)$. Thus the desired count is at most
\[
 \#\{m\le x/d:v_p(C(m))\le v_p(e)\}
 \ll_{d,e}x^{1-\eta_p}(1+\log x)^{v_p(e)}
 \ll_{d,e}x^{1-\eta},
\]
where $\eta_p>0$ is as in Lemma~\ref{lem3.1}; we may take $\eta_p=a_d/(L_d\log p)$ and $\eta=\eta_p/2$.
\end{proof}

\section*{Acknowledgements}
The authors thank the referees for their time and consideration.

\section*{Conflict of interest}
The authors declare that they have no competing interests.

\end{document}